\documentclass[MSNbibl,number,citesort,seceqn,dvips]{arxbj}
\usepackage{upgreek}
\font\tencyi=wncyi10

\aid{0}
\volume{18}
\issue{3}
\pubyear{2012}
\firstpage{1042}
\lastpage{1060}
\doi{10.3150/11-BEJ357}

\makeatletter
\setvaluelist{emailmarks}{*,\textdagger,**,\textdaggerdbl}
\newcommand{\Z}{{\mathbf Z}}
\newcommand{\R}{\mathbf{R}}

\renewcommand{\P}{\mathrm{P}}
\newcommand {\E}{\mathrm{E}}

\newcommand{\1}{\mathbf{1}}
\renewcommand{\d}{\mathrm{d}}

\newcommand{\sG}{\mathcal{G}}
\newcommand{\sL}{\mathcal{L}}
\newcommand{\sP}{\mathcal{P}}
\newcommand{\sI}{\mathcal{I}}
\newcommand{\sA}{\mathcal{A}}
\newcommand{\sF}{\mathcal{F}}
\newcommand{\lip}{\operatorname{Lip}}
\renewcommand{\Re}{\operatorname{Re} }

\newtheorem{stat}{Statement}[section]
\newtheorem{proposition}[stat]{Proposition}
\newtheorem{corollary}[stat]{Corollary}
\newtheorem{theorem}[stat]{Theorem}
\newtheorem{lemma}[stat]{Lemma}
\newproclaim{definition}[stat]{Definition}
\newremark{remark}[stat]{Remark}
\newproclaim{assumptions}[stat]{Assumptions}

\renewcommand{\epsilon}{\varepsilon}
\renewcommand{\bm}{\boldsymbol}
\newcommand{\eqref}[1]{(\ref{#1})}
\makeatother

\begin{document}
\begin{frontmatter}

\title{An asymptotic theory for randomly forced discrete nonlinear heat equations}
\runtitle{The discrete nonlinear heat equation}

\begin{aug}
%%%% inicialai - be tarpu
\author[a]{\fnms{Mohammud} \snm{Foondun}\corref{}\thanksref{e1,u1}\ead[label=e1,mark]{mohammud@math.utah.edu}\ead[label=u1,url,mark]{http://www.math.utah.edu/\textasciitilde mohammud}} \and
\author[a]{\fnms{Davar} \snm{Khoshnevisan}\thanksref{e2,u2}\ead[label=e2,mark]{davar@math.utah.edu}\ead[label=u2,url,mark]{http://www.math.utah.edu/\textasciitilde davar}}
%%\runauthor{} %% auto
\address[a]{Department of Mathematics, University of Utah,
                Salt Lake City, UT 84112-0090, USA.\\\printead{e1,e2}\\\printead{u1,u2}}
\end{aug}

% HISTORY:
\received{\smonth{8} \syear{2010}}
\revised{\smonth{2} \syear{2011}}

% ABSTRACT
\begin{abstract}
We study discrete nonlinear
parabolic stochastic heat equations of the form, $u_{n+1}(x) - u_n(x) =
(\sL u_n)(x) + \sigma(u_n(x))\xi_n(x)$, for $n\in \Z_+$ and
$x\in \Z^d$, where $\bm\xi:=\{\xi_n(x)\}_{n\ge 0,x\in\Z^d}$
denotes random forcing and $\sL$ the
generator of a random walk on $\Z^d$. Under mild
conditions, we  prove that
the preceding stochastic PDE has a unique solution
that grows at most exponentially in time.
And that, under natural conditions,  it  is ``weakly intermittent.''
Along the way, we establish a~comparison principle as well as a~finite support property.
\end{abstract}

% KEYWORDS
\begin{keyword}
\kwd{intermittency}
\kwd{stochastic heat equations}
\end{keyword}

\end{frontmatter}

%s1 #&#
%s1 ###
\section{Introduction}\label{sec1}

Let us consider a prototypical stochastic heat equation
of the following type:
%e1.1 #&#
%e1.1 ###
\begin{equation}\label{heatcont}
\everymath{\displaystyle }
\left|\begin{array}{l}
        \frac{\partial u(t ,x)}{\partial t} = (\sL u)(t ,x) +
                \sigma(u(t ,x))\xi_t(x)
                  \qquad \mbox{for $t>0$ and $x\in\R$},\\\noalign{\vspace*{4pt}}
        u(0 ,x)= u_0(x),
\end{array}\right.\end{equation}
where $u_0$ and $\sigma$ are known non-random functions:
$u_0$ is  bounded and measurable;
\mbox{$\sigma\dvtx \R\to\R$} is Lipschitz continuous;
$\bm{\xi}:=\{\xi_t\}_{t\ge 0}$ is an infinite-dimensional white
noise; and $\sL$ is an operator acting on the variable $x$. It
is well known that \eqref{heatcont}
has a unique ``mild solution''
under natural conditions on  $\bm\xi$ and $\sL$
\cite{Dalang,DZ,Funaki,KR79a,KR79b,KR77,Pardoux75,Pardoux72,Walsh};
we can think of $\bm \xi$ as the ``forcing term'' as well as the ``noise.''

Let us observe that, in \eqref{heatcont},
the operator $\sL$ and the noise term compete with one another:
$\sL$ tends to flatten/smooth the solution $\mathbf u$,
whereas the noise term tends to make~$\mathbf u$ more irregular.
This competition was studied in \cite{FKN}
in the case that $\sigma:=1$ and $\sL:=$ the $L^2$-generator of a L\'evy process.

The [parabolic] ``Anderson model''
is an important special case of \eqref{heatcont}.
In that case one considers $\sL:=\kappa \partial_{xx}$ and $\sigma(z):=\nu z$ for fixed $\nu,\kappa>0$,
and interprets $u(t ,x)$ as the average number of particles -- at
site $x$ and time $t$ -- when the particles perform independent
Brownian motions; every
particle splits into two at rate $\xi_t(x)$ -- when $\xi_t(x)>0$ -- and
is extinguished at rate $-\xi_t(x)$ -- when $\xi_t(x)<0$.
See Carmona and Molchanov \cite{CM}, Chapter~1, for this,
together with a groundbreaking analysis of the ensuing model.
The Anderson model also has important
connections to stochastic analysis, statistical physics,
random media, cosmology, etc.\
\cite{BC,BCJ,CKM,CV98,CM94,CM,CM1,CMS,CMS1,FV,Funaki,GartnerdenHollander,GartnerKonig,GKM,GK,GN,HKM,Kardar,KPZ,KLMS,KrugSpohn,LL63,Molch91,Mueller,SZ,Shiga,ZMRS85,ZMRS88}.  We note that
many of the mentioned papers are concerned with time-independent noise
only.\looseness=-1

A majority of the sizable literature on the Anderson model is concerned
with establishing a property called ``intermittency''
\cite{Mandelbrot,Molch91,Novikov,ZMRS88,ZMRS85}.
Recall that the \emph{pth moment Lyapounov exponent $\gamma(p)$} is
defined as
%e1.2 #&#
%e1.2 ###
\begin{equation}\label{eqgamma}
        \gamma(p):=
        \lim_{t\to\infty} \frac 1t \ln \E  [
        u(t ,x)^p ],
\end{equation}
provided that the limit exists.
The solution $\mathbf{u}:=\{u(t ,x)\}_{t\geq 0, x\in \R^d}$
to the parabolic Anderson model is
said to be \emph{intermittent} if $\gamma(p)$ exists for all $p\ge 1$ and
$p\mapsto (\gamma(p)/p)$ is strictly increasing on $[1 ,\infty)$.
This mathematical definition describes a ``separation of scales'' phenomena
and is believed to capture many of the salient features
of its physical counterpart in statistical physics and turbulence
\cite{BatchelorTownsend,Mandelbrot,Novikov,vonWeizsacker,ZMRS85}.
For more information see the Introductions of Bertini and Cancrini \cite{BC}
and Carmona and Molchanov \cite{CM}.

Recently \cite{FK} we considered \eqref{heatcont}
in a fully nonlinear setting with space--time white noise $\bm\xi$
and $\sL:=$ the $L^2$-generator of a L\'evy process.
We showed that if $\sigma$ is ``asymptotically linear''
and $u_0$ is ``sufficiently large,'' then
$p\mapsto \widetilde{\gamma}(p)/p$ is strictly
increasing on $[2 ,\infty)$, where
%e1.3 #&#
%e1.3 ###
\begin{equation}\label{eqgamma}
        \widetilde{\gamma}(p):=
        \limsup_{t\to\infty} \frac 1t \ln \E  ( |
        u_t(x)|^p ).
\end{equation}
This gives evidence of
intermittency for solutions of stochastic PDEs. Moreover,
bounds on $\widetilde{\gamma}$ were given in
terms of the Lipschitz constant of $\sigma$ and the function
%e1.4 #&#
%e1.4 ###
\begin{equation}\label{defcontUpsilon}
        \widetilde{\Upsilon}(\beta)
        := \frac{1}{2\uppi}
        \int_{-\infty}^\infty \frac{\mathrm{d}\xi}{\beta+2\Re\Psi(\xi)},
         \qquad  \mbox{defined for all $\beta>0$},
\end{equation}
where $\Psi$ denotes the characteristic
exponent of the L\'evy process generated by the $\sL$.
It is precisely this connection between $\widetilde{\Upsilon}$
and $\sigma$ that allows us to
describe a relationship between the smoothing effects
of $\sL$ and the roughening effect of the underlying
forcing terms.

There are two physically relevant classes of bounded initial data
$u_0$ that arise naturally in the literature \cite{BC,CM94,Molch91}:
{(a)}~Where~$u_0$ is bounded below, away from zero; and
{(b)}~Where~$u_0$ has compact support.
Our earlier analysis \cite{FK} studies
fairly completely Case~{(a)} but fails to say anything about Case
{(b)}. We do not know much about {(b)}, in fact.
Our present goal is to consider instead
a \emph{discrete} setting in which we
\emph{are} able to analyze Case~{(b)}.

There is a large literature on [discrete] partial difference
equations of the heat type; see Agarwal \cite{Agarwal} and its many
chapter bibliographies. Except for the work by Zeldovich \textit{et al.}
\cite{ZMRS85}, Section~5, we have found little
on fully-discrete stochastic heat equations~\eqref{heatcont}.
We will see soon that the discrete setup treated here yields many of
the interesting mathematical features that one might wish for,
and at low technical cost.
For instance, we do not presuppose a knowledge of PDEs
and/or stochastic calculus in this paper.

An outline of the paper follows: In Section~\ref{sec2},
we state the main results of the paper; they
are proved in Section~\ref{sec5}, after we establish some auxiliary
results in Section~\ref{sec3} and Section~\ref{sec4}.
In Section~\ref{sec6} we compute a version of the second-moment
[upper] \emph{Lyapounov exponent} of the solution $u$ to the
parabolic Anderson model with temporal noise.
From a physics point of view, that model is only modestly interesting,
but it provides a setting in which
we can rigorously verify many of the predictions of the replica
method \cite{Kardar}. The replica method
itself will not be used, however.

Throughout the paper, we define
%e1.5 #&#
%e1.5 ###
\begin{equation}\label{norm}
        \|X\|_p:=\{\E (|X|^p)\}^{1/p}  \qquad\mbox{for all $X\in L^p(\P)$},
\end{equation}
for every $p\in[1 ,\infty)$.

%s2 #&#
%s2 ###
\section{Main results}\label{sec2}

Throughout, we study
the following discrete version of \eqref{heatcont}:
%e2.1 #&#
%e2.1 ###
\begin{equation}\label{heat}
        u_{n+1}(x) - u_n(x) = (\sL u_n)(x) + \sigma(u_n(x))
        \xi_n(x)\qquad\mbox{for $n\geq0$ and $x\in\Z^d$},
\end{equation}
with [known] bounded initial function $u_0\dvtx \Z^d\to\R$ and diffusion coefficient
$\sigma\dvtx \R\to\R$. The operator
$\sL$ acts on functions of $x$ and is
the generator of a random walk on $\Z^d$.

Let $\sI$ denote the identity operator and
$\sP:=\sL+\sI$ the transition operator for $\sL$.
Then~\eqref{heat} is equivalent to
the following recursive relation:
%e2.2 #&#
%e2.2 ###
\begin{equation}\label{heat1}
        u_{n+1}(x) = (\sP u_n)(x) + \sigma(u_n(x))\xi_n(x).
\end{equation}

Our first contribution is an analysis of \eqref{heat} in the
case that the $\xi$'s are i.i.d. with
common mean 0 and variance 1 [discrete white noise].
The following function \mbox{$\Upsilon\dvtx (1 ,\infty)\to\R_+$}
is the present analogue of $\widetilde{\Upsilon}$
[see \eqref{defcontUpsilon}]:
%e2.3 #&#
%e2.3 ###
\begin{equation}\label{defUpsilon}
        \Upsilon(\lambda) := \frac{1}{(2\uppi)^d}
        \int_{(-\uppi,\uppi)^d} \frac{\d\xi}{\lambda -
        |\phi(\xi)|^2}\qquad
        \mbox{for all $\lambda >1$},
\end{equation}
where $\phi$ denotes the characteristic exponent of
the increments of the walk that corresponds to $\sL$;
that is,
%e2.4 #&#
%e2.4 ###
\begin{equation}
        \phi(\xi) := \sum_{x\in\Z^d} \mathrm{e}^{\mathrm{i}x\cdot \xi} P_{0,x},
\end{equation}
and $P_{0,\bullet}$ is the transition function of the random walk.
Because $\Upsilon$ is continuous,
strictly positive, and strictly
decreasing on $(1 ,\infty)$, it has a continuous
strictly decreasing inverse on $(0 ,\Upsilon(1^-))$.
We extend the definition of that inverse by setting
%e2.5 #&#
%e2.5 ###
\begin{equation}
        \Upsilon^{-1}(x) :=
        \sup \{ \lambda>1\dvtx
        \Upsilon(\lambda)>x
         \},
\end{equation}
where $\sup\varnothing:=1$.
Also, let
%e2.6 #&#
%e2.6 ###
\begin{equation}
        \lip_\sigma:= \sup_{x\neq y}
        \frac{ |\sigma(x)-\sigma(y)|}{|x-y|}
\end{equation}
denote the Lipschitz constant of the function $\sigma$
[$\lip_\sigma$ can be infinite].
The following is a~discrete counterpart
of Theorems 2.1 and 2.7 of \cite{FK}, and is our first main result.

%th2.1 #&#
\begin{theorem}\label{thspacetime}
        Suppose $\bm{\xi}$ are i.i.d. with mean 0 and
        variance 1. If $u_0$ is bounded and~$\sigma$
        is Lipschitz continuous, then \eqref{heat} has an
        a.s.-unique solution $\mathbf{u}$ which satisfies the following:
        For all $p\in[2 ,\infty)$,
        %e2.7 #&#
%e2.7 ###
\begin{equation}\label{eqspacetime}
                \limsup_{n\to\infty} \frac1n \sup_{x\in\Z^d}\ln
                 \| u_n(x)  \|_p
                \le \frac1 2\ln  \Upsilon^{-1}
                 ((c_p\lip_\sigma)^{-2} ),
        \end{equation}
        where $\Upsilon^{-1}(1/0):=0$ and
        $c_p:=$ the optimal constant in Burkholder's
        inequality for discrete-parameter martingales.
        Conversely, if $\inf_{x\in\Z^d}u_0(x)>0$ and
        $L_\sigma:=\inf_{z\in\R}
        |\sigma(z)/z|>0$, then
        %e2.8 #&#
%e2.8 ###
\begin{equation}\label{eqspacetime1}
                \inf_{x\in\Z^d}\limsup_{n\to\infty}  \frac1n \ln
                 \| u_n(x)  \|_p
                \ge \frac12 \ln \Upsilon^{-1}
                 (L_\sigma^{-2} ).
        \end{equation}
\end{theorem}

To avoid confusion, we emphasize that the norm
in the above inequalities is the $L^p(\P)$ norm
that was defined in \eqref{norm}.

The exact value of $c_p$ is not known \cite{Wang}.
Burkholder's method itself produces
$c_p \le 18pq^{1/2}$ \cite{HH}, Theorem 2.10, page 23, where
$q:=p/(p-1)$ denotes the conjugate to $p$. It is likely
that better bounds are known, but we are not aware of them.

Theorem \ref{thspacetime} has a continuous
counterpart in \cite{FK}.  Next we point out that
``$u_0$ is bounded below'' in Theorem \ref{thspacetime}
can sometimes be replaced by ``$u_0$ has finite support.''
As far as we know, this does not seem to have a continuous analogue
\cite{FK}. But first we recall the following standard definition:

%de2.2 #&#
\begin{definition}\label{deflocal}
        $\sL$ is \emph{local} if there exists $R>0$ such that $P_{0,x}=0$
        if $|x|>R$.\footnote{As
        is sometimes customary, we identify the Fredholm
        operator $\sP^n$ with its kernel, which is merely the $n$-step
        transition probability:
        $P^n_{x,y}=P^n_{0,y-x}$ at $(x ,y)\in\Z^d\times\Z^d$.}
\end{definition}

%th2.3 #&#
\begin{theorem}\label{thspacetimebis}
        Suppose $\sL$ is local, and the
        $\xi$s are i.i.d. with mean 0 and
        variance 1. In addition,  $u_0\not\equiv 0$
        has finite support,
        $\sigma$ is Lipschitz continuous with $\sigma(0)=0$, and
        $L_\sigma:=\inf_{z\in\R}
        |\sigma(z)/z|>0$. Then,
        for all $p\in[2 ,\infty)$,
        $M_n:=\sup_{x\in\Z^d}|u_n(x)|$ satisfies
        %e2.9 #&#
%e2.9 ###
\begin{eqnarray}
        \label{eqspacetimebis}
                 \frac12\ln \Upsilon^{-1}
                         ( L_\sigma^{-2} )
                        &\le& \limsup_{n\to\infty} \frac1n
                        \sup_{x\in\Z^d}\ln\|u_n(x)\|_p\nonumber
                        \\[-8pt]
                        \\[-8pt]
                & \le& \limsup_{n\to\infty} \frac1n \ln
                         \| M_n  \|_p
                        \le \frac12 \ln \Upsilon^{-1}
                         ((c_p\lip_\sigma)^{-2} ).
        \nonumber
\end{eqnarray}
\end{theorem}
%
%       When $\sigma(z):=\nu z$ [Anderson model],
%       Theorem \ref{thspacetime} yields a closed-form formula
%       for the ``Lyapounov exponent''
%       $\limsup_{n\to\infty} n^{-1}\sup_x\ln\|u_n(x)\|_{L^2(\P)}$.
%       It turns out that this constant can be hard to understand
%       in the simplest settings; see Proposition \ref{prU} for details.
%       \qed

Under the assumption that $\mathbf{u}:=\{u_n(x) \}_{n>0,x\in \Z^d}$ is non-negative,
we define the \emph{upper pth-moment Lyapounov exponent}
$\bar{\gamma}(p)$ as follows:
%e2.10 #&#
%e2.10 ###
\begin{equation}\label{defliap}
        \bar{\gamma}(p):=\limsup_{n \rightarrow \infty} \frac1n
        \sup_{x\in \Z^d} \ln \E [ u_n(x)^p ].\vspace*{-2pt}
\end{equation}

%de2.4 #&#
\begin{definition}
        We say that $\mathbf{u}:=\{u_n(x) \}_{n>0,x\in \Z^d}$ is
        \emph{weakly intermittent} if
        $\bar{\gamma}(p)<\infty$ for all positive and finite $p$,
        and $p\mapsto(\bar\gamma(p)/p)$ is strictly increasing on
        $[2 ,\infty)$.\vspace*{-2pt}
\end{definition}

Our next result is a consequence of the previous theorem, and assumes -- among other things -- that $\sigma(0)=0$.
This condition ensures that the solution has finite support.\vspace*{-2pt}
%co2.5 #&#
\begin{corollary}\label{corinter}
        Suppose, in addition to the conditions of Theorem \ref{thspacetimebis},
        that
\[
        C_\xi:=\sup_{n\geq 0}\sup_{x\in \Z^d}|\xi_n(x)|
        \]
         is finite and
        %e2.11 #&#
%e2.11 ###
\begin{equation}\label{eqcomparison0}
                P_{0,0}\ge C_\xi\lip_\sigma.
        \end{equation}
        Then \eqref{heat} has a weakly intermittent solution.\vspace*{-2pt}
\end{corollary}

We emphasize that $\bar{\gamma}(p)$ is
\emph{not} an exact discrete version of $\widetilde{\gamma}(p)$,
as it is missing absolute values. Condition \eqref{eqcomparison0} allows $\bar{\gamma}(p)$
to be well defined for our solution.

Our next result concerns the Anderson model with temporal noise.
In other words we consider \eqref{heat} with $\sigma(z)=z$,
$\xi_n(x)=\xi_n$ for all $x\in\Z^d$, and
$\bm{\xi}:=\{\xi_n\}_{n=0}^\infty=$
i.i.d. random variables. The present model is
motivated by Mandelbrot's analysis of random
cascades in turbulence \cite{Mandelbrot} and is designed to
showcase a family of examples where the predictions of
the replica method of Kardar \cite{Kardar} can be shown rigorously.
We make the following assumptions:\vspace*{-2pt}

%as2.6 #&#
\begin{assumptions}\label{condtemporal}
        Suppose:
        \begin{enumerate}[(b)]
                \item[(a)] $\sL$ is \emph{local};
                \item[(b)] $\sup_{n\ge 0}|\xi_n| \le P_{0,0}<1$
                [lazy, non-degenerate random walk]; and
        \item[(c)] $u_0(x)\ge 0$ for all $x\in\Z^d$,
                and $0<\sum_{z\in\Z^d} u_0(z)<\infty$.\vspace*{-2pt}
        \end{enumerate}
\end{assumptions}

Then we offer the following.\vspace*{-2pt}
%
%th2.7 #&#
\begin{theorem}\label{thtemporal}
        Under Assumptions \ref{condtemporal}, the Anderson model
        $[$\eqref{heat} with $\sigma(z)=z]$ has a unique
        a.s.-nonnegative solution $\mathbf{u}$, and
        $ M_n:=\sup_{x\in\Z^d} u_n(x)$ satisfies
     %e2.12 #&#
%e2.12 ###
\begin{eqnarray}
                \lim_{n\to\infty} \frac1n \ln\E (
                        M_n^p ) &=&\Gamma(p)
                        \qquad\mbox{for all $p\in[0 ,\infty)$} \quad
                        \mbox{and}\nonumber
                        \\[-8pt]
                        \\[-8pt]
                        \lim_{n\to\infty} \frac1n \ln  M_n
                &=&\lim_{n\to\infty} \frac1n \E (\ln  M_n
                         ) = \Gamma'(0^+)
                        \qquad\mbox{almost surely},
    \nonumber
\end{eqnarray}
        where $\Gamma(p) :=\ln\E [ ( 1+\xi_1 )^p ]$
        for all $p\geq0$.\vadjust{\goodbreak}
\end{theorem}

%s3 #&#
%s3 ###
\section{Some auxiliary results}\label{sec3}

Let us start with a simple existence/growth lemma.
Note that $\sigma$ is not assumed
to be Lipschitz continuous, and the $\xi$s
need not be random. The proof is not demanding,
but the result itself is unimprovable (Remark \ref{remgrowth}).

%le3.1 #&#
\begin{lemma}\label{lemfirst}
        Suppose there exist finite
        $C_\sigma$ and $\widetilde{C}_\sigma$
        such that       $|\sigma(z)| \le C_\sigma|z| + \widetilde{C}_\sigma$
        for all $z\in\R$. Suppose also
        that $u_0$ is bounded and
        $C_\xi:=\sup_{n\ge 0}\sup_{x\in\Z^d}|\xi_n(x)|$
        is finite. Then~\eqref{heat} has a unique solution
        $\mathbf u$ that satisfies
        %e3.1 #&#
%e3.1 ###
\begin{equation}\label{eqfirst}
                \limsup_{n\to\infty} \frac1n \sup_{x\in\Z^d} {\ln}
                |u_n(x)| \le \ln(1+C_\sigma C_\xi).
        \end{equation}
\end{lemma}

\begin{pf}
        Clearly, $\|\sP h\|_\infty\le\|h\|_\infty$, where
        $\|h\|_\infty$ denotes the supremum norm of a
        function $h\dvtx \Z^d\to\R$. Consequently, \eqref{heat1} implies
        that
        %e3.2 #&#
%e3.2 ###
\begin{equation}
                \|u_{n+1}\|_\infty \le\|u_n\|_\infty
                 (1+C_\sigma C_\xi  ) + \widetilde{C}_\sigma C_\xi.
        \end{equation}
%       We can iterate the preceding to find that
%       \begin{equation}
%               \|u_{n+1}\|_\infty \le \|u_0\|_\infty
%                ( 1+C_\sigma C_\xi )^{n+1}
%               +\widetilde{C}_\sigma C_\xi\sum_{j=0}^n
%                ( 1+ C_\sigma C_\xi )^j.
%       \end{equation}
%       The lemma follows immediately from this and \eqref{heat1}.
        We iterate this and apply \eqref{heat1} to conclude the proof.
\end{pf}

%re3.2 #&#
\begin{remark}\label{remgrowth}
        Consider \eqref{heat} with
        $u_0(x)\equiv 1$, $\sigma(z)=z$, and $\xi_n(x)\equiv 1$.
        Then $u_n(x)=2^n$ for all $n\ge 0$ and $x\in\Z^d$,
        and \eqref{eqfirst} is manifestly an identity.
        The results of the \hyperref[sec1]{Introduction} show that when the
        $\xi$'s are mean-zero and independent, then the
        worst-case rate in \eqref{eqfirst} can be improved upon;
        this is another evidence of intermittency.
\end{remark}

The following covers the case when $\xi$'s are random variables.
This existence/growth result is proved in the same
manner as Lemma \ref{lemfirst}; we omit the
elementary proof, and also mention that the following cannot be improved
upon.

%le3.3 #&#
\begin{lemma}\label{lemsecond}
        Suppose there exist finite
        $C_\sigma$ and $\widetilde{C}_\sigma$
        such that $|\sigma(z)| \le C_\sigma|z| + \widetilde{C}_\sigma$
        for all $z\in\R$. Suppose also
        that $u_0(x)$ and $\|\xi_n(x)\|_p$ are bounded uniformly
        in $x\in\Z^d$, and $n\ge 0$,
        for some $p\in[1 ,\infty]$.
        Then \eqref{heat} has an a.s.-unique solution
        $\mathbf u$ that satisfies
        %e3.3 #&#
%e3.3 ###
\begin{equation}
                \limsup_{n\to\infty} \frac1n \sup_{x\in\Z^d} {\ln}
                \| u_n(x)\|_p \le \ln
                 ( 1+C_\sigma K_{p,\xi}  ),
        \end{equation}
        where $K_{p,\xi}:=\sup_{n\ge 0}\sup_{x\in\Z^d}
        \|\xi_n(x)\|_p$.
\end{lemma}

%s3.1 #&#
%s3.1 ###
\subsection{A finite-support property}\label{sec3.1}
Next we demonstrate that the solution to \eqref{heat} has
a \emph{finite-support} property.
A remarkable result of Mueller \cite{Mueller} asserts that
Theorem \ref{thfinitespt} below does not have a
naive continuum-limit analogue.
The present work is closer in spirit to the compact-support theorem
of Mueller and Perkins \cite{MuellerPerkins}.

Let us consider the heat equation \eqref{heat}
and suppose that it has a unique
solution $\mathbf{u}:=\{u_n(x)\}_{n\ge 0,x\in\Z^d}$.
We say that a function $f\dvtx \Z^d\to\R$ has
\emph{finite support} if $\{x\in\Z^d\dvtx
f(x)\neq 0\}$ is finite.  Define
%e3.4 #&#
%e3.4 ###
\begin{equation}
        R_n := \inf \{ r>0\dvtx  u_n(x)=0 \mbox{ for all
        $x\in\Z^d$ such that $|x|>r$} \},
\end{equation}
and let $R$ denote the \emph{radius of support of $\sP$}; that is,
%e3.5 #&#
%e3.5 ###
\begin{equation}\label{defR}
        R:= \inf \{r >0\dvtx   P_{0,x}=0
         \mbox{ for all $x\in\Z^d$ with $|x|>r$} \}.
\end{equation}

%th3.4 #&#
\begin{theorem}\label{thfinitespt}
        Suppose $\sigma(0)=0$, and $\sL$ is local. If,
        in addition,
        $u_0$ has finite support, then so does $u_n$
        for all $n\ge 1$. In fact,
        %e3.6 #&#
%e3.6 ###
\begin{equation}\label{eqgrowspt}
                \# \{ x\in\Z^d\dvt   u_{n+1}(x)\neq 0 \}
                \le 2^d[(n+1)R+R_0]^d
                \qquad\mbox{for all $n\ge 0$}.
        \end{equation}
\end{theorem}

\begin{pf}
        Suppose there exists $n\ge 0$ such that
        $u_n(x)=0$ for all but a finite number of
        points $x\in\Z^d$. We propose to prove that
        $u_{n+1}$ enjoys the same finite-support property.
        This clearly suffices to prove the theorem.
        Because $u_n(x)=0$ for all but a finite number of
        $x$s, \eqref{heat1} tells us that for all
        but a finite number of points
        $x\in\Z^d$, $u_{n+1}(x) %= (\sP u_n)(x)
        = \sum_{{y\in\Z^d: |y-x|\le R}} P_{x,y} u_n(y)$.
        Thus if $u_n$ has finite support, then
        so does $u_{n+1}$, and
        $R_{n+1}\le R+R_n$. Equation \eqref{eqgrowspt} also follows from
        this.
\end{pf}

The locality of $\sL$ cannot be dropped altogether. This
general phenomenon appears earlier. For instance,
Iscoe \cite{Iscoe} showed that the super Brownian motion
has a finite-support property,
and Evans and Perkins \cite{EP} proved
that there Iscoe's theorem does not
hold if the underlying motion is non-local.

%s3.2 #&#
%s3.2 ###
\subsection{A comparison principle}\label{sec3.2}
The result of this subsection is a discrete analogue of
Mueller's well-known and deep comparison principle \cite{Mueller};
but the proof uses very simple ideas.
Throughout we assume that there exist unique solutions
$\mathbf{v}$ and $\mathbf{u}$ to \eqref{heat}
with respective initial data $v_0$ and $u_0$ and that $\sigma\dvtx \R\to\R$ is globally
Lipschitz with optimal Lipschitz constant $\lip_\sigma$.

%th3.5 #&#
\begin{theorem}\label{thcomparison}
        Suppose $C_\xi:=\sup_{n\ge 0}\sup_{x\in\Z^d}|\xi_n(x)|$
        is finite and satisfies
        %e3.7 #&#
%e3.7 ###
\begin{equation}\label{eqcomparison}
                P_{0,0}\ge C_\xi\lip_\sigma.
        \end{equation}
        Then $u_0 \ge v_0$ implies that
        $u_n \ge v_n $ for all $n\ge 0$.
\end{theorem}

In the continuous setting, one usually appeals to Mueller's comparison principle,
using a condition such as ``$\sigma(0)=0$,'' in order to establish positivity of
the solution. Therefore, it might be worth
noting that the preceding does \textit{not} require that $\sigma(0)=0$.

\begin{pf*}{Proof of Theorem~\ref{thcomparison}}
        We propose to prove that if $u_n\ge v_n$,
        then  $u_{n+1}\ge v_{n+1}$.
        Let us write $f_k:= u_k-v_k$ for
        all $k\ge 0$. By \eqref{heat1} and \eqref{eqcomparison},
 %e3.8 #&#
%e3.8 ###
\begin{eqnarray}
                f_{n+1}(x)
                &=& (\sP f_n)(x) +  [\sigma(u_n(x))-\sigma(v_n(x)) ]
                        \xi_n(x)\nonumber\\
                &\ge& (\sP f_n)(x) - \lip_\sigma\cdot |f_n(x) \xi_n(x)|\\
                &\ge& (\sP f_n)(x) -  P_{0,0}\cdot |f_n(x)|.\nonumber
    \end{eqnarray}

        But $(\sP h)(x) =\sum_{y\in\Z^d} P_{x,y} h(y)
        \ge P_{0,0} \cdot h(x)$ for all $x\in\Z^d$,
        as long as $h\ge 0$.
        By the induction hypothesis, $f_n$ is
        a non-negative function, and hence so is $f_{n+1}$.
        This gives the desired result.
\end{pf*}

The following ``positivity principle'' follows readily
from Theorem \ref{thcomparison}.

%co3.6 #&#
\begin{corollary}\label{corcomparison}
        If $u_0\geq 0$ in Theorem \ref{thcomparison},
        then $u_n\geq 0$ for all $n\geq 0$.
\end{corollary}

%s4 #&#
%s4 ###
\section{A priori estimates}\label{sec4}
In this section we develop some tools needed for
the proof of Theorems \ref{thspacetime} and \ref{thspacetimebis}.
It might help to emphasize that we are considering
the case where the random field
$\bm{\xi}:=\{\xi_n(x)\}_{n\ge 0,x\in\Z^d}$ is
[discrete] \emph{white noise}. That is,
the $\xi$'s are mutually independent and
have mean 0 and variance 1. [In fact, they will not be assumed to
be identically distributed.] Note, in particular,
that $K_{1,\xi}=0$ and $K_{2,\xi}=1$, where $K_{1,\xi}$
and $K_{2,\xi}$ were defined in Lemma \ref{lemsecond}.

Here, and throughout, let $\sG:=\{\sG_n\}_{n=0}^\infty$
denote the filtration generated by the infinite-dimensional
``white-noise'' $\{\xi_n\}_{n=0}^\infty$. Recall that
a random field $\mathbf{f}:=\{f_n(x)\}_{n\ge 0,x\in\Z^d}$
is \emph{$\sG$-predictable} if the random function
$f_n$ is measurable with respect to $\sG_{n-1}$ for all
$n\ge 1$, and $f_0$ is non-random.

Given a $\sG$-predictable random field $\mathbf{f}$ and
$\lambda>1$, we define
%e4.1 #&#
%e4.1 ###
\begin{equation}\label{defnorm}
        \|\mathbf{f}\|_{\lambda, p} := \sup_{n\ge 0}\sup_{x\in\Z^d}
        \lambda^{-n}\| f_n(x) \|_p.
\end{equation}

Define for all $\sG$-predictable random fields $\mathbf{f}$,
%e4.2 #&#
%e4.2 ###
\begin{equation}\label{eqA}
        (\sA \mathbf{f})_n(x) := \sum_{j=0}^n\sum_{y\in\Z^d}
        P^{n-j}_{x,y}\sigma(f_j(y))\xi_j(y).
\end{equation}

We begin by developing an a priori estimate
for the operator norm of $\sA$.
This estimate
is a discrete $L^p$-counterpart of Lemma
3.3 in \cite{FK}, while the continuity estimates
given by Proposition~\ref{prNA2} is a discrete version of Lemma 3.4 in \cite{FK}.
%pr4.1 #&#
\begin{proposition}\label{prNA1}
        For all $\sG$-predictable random fields $\mathbf{f}$ and all
        $\lambda>1$,
        %e4.3 #&#
%e4.3 ###
\begin{equation}
                \|\sA\mathbf{f}\|_{\lambda,p} \le c_p
                \bigl(|\sigma(0)|+\lip_{\sigma}\|\mathbf{f}\|_{\lambda,p}\bigr)
                \cdot\sqrt{\lambda^2\Upsilon(\lambda^2)}.
        \end{equation}
\end{proposition}

The proof requires a simple arithmetic result
\cite{FK}, Lemma 3.2:

%le4.2 #&#
\begin{lemma}\label{lemarith}
        $(a+b)^2 \le (1+\epsilon)a^2 +( 1+\epsilon^{-1}) b^2$
        for all $a,b\in\R$ and $\epsilon>0$.
\end{lemma}

Hereforth, define
%e4.4 #&#
%e4.4 ###
\begin{equation}\label{eqq}
        q_k := \sum_{z\in\Z^d}  | P^k_{0,z}  |^2
        \qquad\mbox{for all $k\ge 0$}.
\end{equation}
The proof of Proposition \ref{prNA1} also requires the
following Fourier-analytic interpretation of the function
$\Upsilon$.

%le4.3 #&#
\begin{lemma}\label{lemUpsilon}
        $\lambda \Upsilon(\lambda)
        =\sum_{n=0}^\infty \lambda^{-n} q_n$
        for all $\lambda>1$.
\end{lemma}

\begin{pf}
        By the Plancherel theorem \cite{Rudin}, page 26,
        %e4.5 #&#
%e4.5 ###
\begin{equation}
                q_n=\frac{1}{(2\uppi)^d}\int_{(-\uppi,\uppi)^d}
                |\phi(\xi)|^{2n} \,\d\xi.
        \end{equation}
        Multiply the preceding by $\lambda^{-n}$ and
        add over all $n\ge 0$ to finish.
\end{pf}

\begin{pf*}{Proof of Proposition \ref{prNA1}}
        According to Burkholder's inequality,
        %e4.6 #&#
%e4.6 ###
\begin{equation}\label{eqpreBurk}
                \E ( |(\sA\mathbf{f})_n(x)|^p )
                \leq c_p^p\E\Biggl ( \Biggl|
                \sum_{j=0}^n\sum_{y\in \Z^d}|P^{n-j}_{x,y}|^2
                \cdot |\sigma(f_j(y))|^2
                \Biggr |^{p/2} \Biggr).
        \end{equation}
        Since $p/2$ is a positive integer, the
        generalized H\"older inequality  yields
        the following: For all $j=0,\ldots,n$ and
        $y_1,\ldots,y_{p/2}\in\Z^d$,
        %e4.7 #&#
%e4.7 ###
\begin{equation}
                \E  \Biggl(\prod_{i=1}^{p/2} | \sigma(f_j(y_i))|^2  \Biggr)
                \leq \prod_{i=1}^{p/2}  \|\sigma(f_j(y_i))
                 \|_p^2.
        \end{equation}
        After a little algebra, the preceding and \eqref{eqpreBurk}
        together imply that
        %e4.8 #&#
%e4.8 ###
\begin{equation}
                \| (\sA\mathbf{f})_n(x)\|_p^2\leq c_p^2
                \sum_{j=0}^n\sum_{y\in \Z^d}|P^{n-j}_{x,y}|^2\cdot
                \|\sigma(f_j(y))\|_p^2.
        \end{equation}
        Because $\sigma$ is Lipschitz,
        $|\sigma(x)|\le |\sigma(0)|+\lip_\sigma|x|$ for
        all $x\in\R$. Consequently, by~Lem\-ma~\ref{lemarith} and Minkowski's inequality,
             %e4.9 #&#
%e4.9 ###
\begin{eqnarray}
                         \| (\sA\mathbf{f})_n(x) \|^2_p
                                &\le& c_p^2(1+\epsilon)|\sigma(0)|^2
                                \sum_{k=0}^n q_k\nonumber
                                \\[-8pt]
                                \\[-8pt]
                        && {}+c_p^2 (1+\epsilon^{-1} )\lip_\sigma^2
                                \sum_{j=0}^n \sum_{y\in\Z^d}
                                |P^{n-j}_{x,y}|^2\cdot\|f_j(y)\|_p^2.
                \nonumber
 \end{eqnarray}
        In accord with Lemma \ref{lemUpsilon},
        %e4.10 #&#
%e4.10 ###
\begin{equation}
                \sum_{k=0}^nq_k
                \le \lambda^{2n+2} \Upsilon(\lambda^2)
                \qquad\mbox{for all $n\ge 0$}
        \end{equation}and also
        %e4.11 #&#
%e4.11 ###
\begin{equation}\label{eqA3}
                \sup_{y\in\Z^d} \| f_j(y)\|^2_p
                \le \lambda^{2j}\|\mathbf{f}\|_{\lambda,p}^2
                \qquad\mbox{for all $j\ge 0$}.
        \end{equation}
        It follows that $\lambda^{-2n} \| (\sA\mathbf{f})_n(x) \|^2_p$
        is bounded above by
%e4.12 #&#
%e4.12 ###
\begin{eqnarray}
                &&c_p^2(1+\epsilon)|\sigma(0)|^2\lambda^2\Upsilon(\lambda^2)\nonumber
                \\[-8pt]
                \\[-8pt]
                && \quad  {}+ c_p^2 (1+\epsilon^{-1} )\lip_\sigma^2
                        \sum_{j=0}^n \sum_{y\in\Z^d}
                        |P^{n-j}_{x,y}|^2\lambda^{-2(n-j)}\|\mathbf{f}\|_{\lambda,p}^2.
       \nonumber
\end{eqnarray}
        We now take supremum over all $n\geq 1$ and $x\in \Z^d$, and obtain
        %e4.13 #&#
%e4.13 ###
\begin{equation}
                \| \sA\mathbf{f}\|_{\lambda, p}^2\leq c_p^2
                \lambda^2\Upsilon(\lambda^2)\cdot \{
                (1+\epsilon)|\sigma(0)|^2
                + ( 1+\epsilon^{-1} ) \lip^2_{\sigma}\|\mathbf{f}\|_{\lambda,p}^2
                 \}.
        \end{equation}
        We obtain the result upon optimizing the right-hand side over $\epsilon>0$.
\end{pf*}

Next we present an a priori estimate of
the degree to which $\sA$ is continuous.

%pr4.4 #&#
\begin{proposition}\label{prNA2}
        For all predictable random fields $\mathbf{f}$ and $\mathbf{g}$,
        and all $\lambda>1$,
        %e4.14 #&#
%e4.14 ###
\begin{equation}
                \|\sA\mathbf{f} - \sA \mathbf{g}\|_{\lambda, p} \le c_p
                {\lip_\sigma} \|\mathbf{f}-\mathbf{g}\|_{\lambda, p}\cdot\sqrt{%
                \lambda^2\Upsilon(\lambda^2)}.
        \end{equation}
\end{proposition}

\begin{pf}
        We can, and will, assume without loss of generality
        that $\|\mathbf{f}-\mathbf{g}\|_{\lambda,p}<\infty$; else,
        there is nothing to prove.
        By using Burkholder's inequality and arguing as in the previous lemma, we
        find that
        %e4.15 #&#
%e4.15 ###
\begin{equation}
                \|(\sA\mathbf{f})_n(x)-(\sA \mathbf{g})_n(x) \|_p^2
                \le c_p^2\lip_\sigma^2\cdot
                \sum_{j=0}^n\sum_{y\in \Z^d} |P_{x,y}^{n-j}|^2
                \cdot\|f_j(y)-g_j(y) \|_p^2.
        \end{equation}
        We can apply \eqref{eqA3}, but with $\mathbf{f}-\mathbf{g}$
        in place of $\mathbf f$,
        and follow the proof of  Lemma
        \ref{lemUpsilon} to finish the proof.
\end{pf}

%s5 #&#
%s5 ###
\section{Proof of main results}\label{sec5}
Before we prove the main results
we provide a version of Duhamel's principle for discrete equations.
%pr5.1 #&#
\begin{proposition}[(Duhamel's principle)]\label{prDuhamel}
        Suppose that there exists a unique solution to \eqref{heat};
        then for all $n\ge 0$ and $x\in\Z^d$,
        %e5.1 #&#
%e5.1 ###
\begin{equation}
                u_{n+1}(x) = (\sP^{n+1}u_0)(x) + \sum_{j=0}^n
                \sum_{y\in\Z^d} P^{n-j}_{x,y} \sigma(u_j(y))\xi_j(y)
                  \qquad \mbox{a.s.}
        \end{equation}
\end{proposition}

%re5.2 #&#
\begin{remark}
        Among other things, Proposition \ref{prDuhamel} implies
        that $\{u_n\}_{n=0}^\infty$ is an infinite-dimensional Markov chain
        with values in $(\Z^d)^{\Z_+}$ and
        that  $u_{n+1}$ is measurable with respect to
        $\{\xi_k(\bullet)\}_{k=0}^n$ for all $n\ge 0$.
\end{remark}

\begin{pf*}{Proof of Proposition~\ref{prDuhamel}}
        One checks directly that \eqref{heat1} implies that
        $(\sP u_n)(x)$ can be written as
        $(\sP^2 u_{n-1})(x)+\sum_{y\in\Z^d} P_{x,y}
        \sigma(u_{n-1}(y))\xi_{n-1}(y)$,
        and the proposition follows a simple induction
        scheme.
\end{pf*}

%s5.1 #&#
%s5.1 ###
\subsection{Remaining proofs}\label{sec5.1}
\begin{pf*}{Proof of Theorem \ref{thspacetime}}
        We proceed in two steps:
        First we prove uniqueness and \eqref{eqspacetime},
        and then we establish \eqref{eqspacetime1}.

        \emph{Step 1:} Let
        $f^{(0)}_n(x) := u_0(x)$ for all $n\ge 0$ and $x\in\Z^d$.
        We recall the operator $\sA$ from~\eqref{eqA},
        and define iteratively a predictable
        random field $\mathbf{f}^{(\ell+1)}$
        as follows: $f^{(\ell+1)}_0(x):=u_0(x)$ for all
        $x\in\Z^d$, and
        %e5.2 #&#
%e5.2 ###
\begin{equation}\label{deffell}
                f^{(\ell+1)}_{n+1}(x) := (\sP^{n+1}u_0)(x)
                +  \bigl(\sA\mathbf{f}^{(\ell)} \bigr)_n(x),
        \end{equation}
        for integers $n,\ell\ge 0$ and $x\in\Z^d$.
        Proposition \ref{prNA1} and induction together imply
        that $\|\sA\mathbf{f}^{(\ell)}\|_{\lambda,p}<\infty$ for all
        $\lambda>1$ and $\ell\ge 0$. And therefore,
        $\|\mathbf{f}^{(m)}\|_{\lambda,p}<\infty$ for all $m\ge 0$
        and $\lambda>1$, as well.
        We multiply \eqref{deffell} by $\lambda^{-n}$ and
        use the fact that $f_0^{(m)}\equiv u_0$ to obtain
        %e5.3 #&#
%e5.3 ###
\begin{equation}
                \bigl\|\mathbf{f}^{(\ell+1)}-\mathbf{f}^{(\ell)} \bigr\|_{\lambda,p}=
                \frac{1}{\lambda}
                \bigl\|\sA\mathbf{f}^{(\ell)}-\sA\mathbf{f}^{(\ell-1)}\bigr\|_{\lambda,p}.
        \end{equation}
        Thus, Proposition \ref{prNA2} implies
        %e5.4 #&#
%e5.4 ###
\begin{equation}
                \bigl\|\mathbf{f}^{(\ell+1)}-\mathbf{f}^{(\ell)} \bigr\|_{\lambda,p}\leq c_p
                \lip_{\sigma}\sqrt{\Upsilon(\lambda^2)}
                \cdot\bigl\|\mathbf{f}^{(\ell)}-
                \mathbf{f}^{(\ell-1)} \bigr\|_{\lambda,p}.
        \end{equation}
        This and iteration together yield
        %e5.5 #&#
%e5.5 ###
\begin{equation}\label{eqf1,0}
                \bigl\|\mathbf{f}^{(\ell+1)}-\mathbf{f}^{(\ell)} \bigr\|_{\lambda,p}\leq
                 \bigl( c_p\lip_{\sigma}\sqrt{\Upsilon(\lambda^2)}
                 \bigr)^\ell \cdot \bigl\| \mathbf{f}^{(1)}-\mathbf{f}^{(0)} \bigr\|_{\lambda,p}.
        \end{equation}
        In order to estimate the final $(\lambda ,p)$-norm we use
        \eqref{deffell} $[\ell:=0]$
        and Minkowski's inequality to find that
        %e5.6 #&#
%e5.6 ###
\begin{equation}
                \bigl\|f_{n+1}^{(1)}-f_{n+1}^{(0)} \bigr\|_p
                \le 2\|u_0(x)\|_p+\bigl\|\bigl(\sA \mathbf{f}^{(0)}\bigr)_n \bigr\|_p.
        \end{equation}
        We argue as before and use Proposition \ref{prNA1} to deduce
        that  $\|\mathbf{f}^{(1)}-\mathbf{f}^{(0)}\|_{\lambda, p}$
        is bounded above by $2\|u_0\|_{\lambda,p}
        +c_p(|\sigma(0)|+\lip_{\sigma}\|u_0\|_{\lambda,p})
        \sqrt{\Upsilon(\lambda^2)}$.
        Thus, by \eqref{eqf1,0},
        %e5.7 #&#
%e5.7 ###
\begin{eqnarray}
                &&\bigl\| \mathbf{f}^{(\ell+1)}-\mathbf{f}^{(\ell)}\bigr\|_{\lambda,p}
                \nonumber
                \\[-8pt]
                \\[-8pt]
                && \quad  \le  \bigl(c_p\lip_\sigma\sqrt{\Upsilon(\lambda^2)}
                         \bigr)^\ell\cdot
                         \bigl\{c_p\bigl(|\sigma(0)|+\lip_{\sigma}\|u_0\|_{\lambda,p}\bigr)
                        \sqrt{\Upsilon(\lambda^2)}+2\|u_0\|_{\lambda, p} \bigr\}.
        \nonumber
        \end{eqnarray}
        Consequently, if $\Upsilon(\lambda^2)<(c_p\lip_\sigma)^{-2}$,
        then $\|f^{(\ell+1)}-f^{(\ell)}\|_{\lambda,p}$ is summable in
        $\ell$. Whence there exists a predictable $\mathbf{f}$ such that
        $\|\mathbf{f}^{(\ell)}-\mathbf{f}\|_{\lambda,p}$ tends to zero as $\ell$ tends to
        infinity, and~$\mathbf{f}$ solves \eqref{heat}.
        Proposition \ref{prDuhamel} implies that
        $f_n(x)=u_n(x)$ a.s.,
        for all $n\ge 0$ and $x\in\Z^d$.
        It follows that $\|\mathbf{u}\|_{\lambda,p}<\infty$
        provided that $\Upsilon(\lambda^2)<(c_p\lip_\sigma)^{-2}$.
        The first part of the theorem -- that is, existence,
        uniqueness, and \eqref{eqspacetime} -- all
        follow from this finding. We now turn our attention to the second step of the proof.

        \emph{Step 2:}
        Hereforth, we assume that
        $\alpha:=\inf_{x\in\Z^d} u_0(x)>0$
        and $|\sigma(z)|\ge L_\sigma|z|$ for all $z\in\R$.
        It follows readily from Proposition \ref{prDuhamel} that
    %e5.8 #&#
%e5.8 ###
\begin{eqnarray}\label{eqll1}
                \E (  | u_{n+1}(x) |^2 ) &=&
                         | (\sP^{n+1}u_0)(x) |^2 + \sum_{j=0}^n
                        \sum_{y\in\Z^d}  | P^{n-j}_{x,y} |^2
                        \E ( | \sigma(u_j(y)  |^2 )\nonumber
                        \\[-8pt]
                        \\[-8pt]
        %       &\ge \alpha^2 + L_\sigma^2 \cdot\sum_{j=0}^n\sum_{y\in\Z^d}
        %                | P^{n-j}_{x,y}  |^2 \E (  |
        %               u_j(y)  |^2 )\\
        %       &= \alpha^2 + L_\sigma^2 \cdot\sum_{j=0}^n\sum_{y\in\Z^d}
        %                | P^{n-j}_{0,y-x}  |^2 \E (  |
        %               u_j(y)  |^2 ).
                &\ge& \alpha^2 + L_\sigma^2 \cdot\sum_{j=0}^n\sum_{y\in\Z^d}
                         | P^{n-j}_{0,y-x}  |^2 \E (  |
                        u_j(y)  |^2 ).
        \nonumber
\end{eqnarray}
        In order to solve this we define for all $\lambda>1$ and $z\in\Z^d$,
        %e5.9 #&#
%e5.9 ###
\begin{equation}
                \sF_\lambda(z) := \sum_{j=0}^\infty\lambda^{-j} |
                P^j_{0,z}|^2
                \quad\mbox{and}\quad
                \sG_\lambda(z) := \sum_{j=0}^\infty \lambda^{-j} \E  ( |
                u_j(z) |^2  ).
        \end{equation}
        We can multiply the extreme quantities in \eqref{eqll1}
        by $\lambda^{-(n+1)}$ and add $[n\ge 0]$ to find that
        %e5.10 #&#
%e5.10 ###
\begin{equation}\label{eqll2}
                \sG_\lambda(x) \ge \frac{\alpha^2 \lambda}{\lambda-1} +
                \frac{L_\sigma^2}{\lambda} \cdot  (
                \sF_\lambda* \sG_\lambda )(x),
        \end{equation}
        where the left-hand side is obtained after adding the initial term.
        This is a renewal inequality  \cite{CD}; we prove that \eqref{eqll2} does
        not have a finite solution when $\Upsilon(\lambda)\ge L_\sigma^{-2}$.
        If $\1(x):=1$ for all $x\in\Z^d$, then
        $(\sF_\lambda*\1)(x) = \lambda\Upsilon(\lambda)$
        for all $x\in\Z^d$ [Lemma \ref{lemUpsilon}].
        Therefore,~\eqref{eqll2} yields
%e5.11 #&#
%e5.11 ###
\begin{eqnarray}
                \sG_\lambda(x)
                        &\ge& \frac{\alpha^2 \lambda}{\lambda -1} +
                        \frac{L_\sigma^2}{\lambda}\cdot\biggl (
                        \frac{\alpha^2 \lambda}{\lambda-1}
                        (\sF_\lambda*\1)(x)+ \frac{L_\sigma^2}{\lambda}
                        \cdot(\sF_\lambda*\sF_\lambda*\sG_\lambda)(x)
                         \biggr)\nonumber
                         \\[-8pt]
                         \\[-8pt]
                &=& \frac{\alpha^2 \lambda}{\lambda-1}  \{ 1+ \Upsilon(\lambda)
                        L_\sigma^2  \} +   \biggl(\frac{L_\sigma^2}{\lambda} \biggr)^2
                        \cdot (\sF_\lambda*\sF_\lambda*\sG_\lambda)(x).
         \nonumber
\end{eqnarray}
        We now use \eqref{eqll2} and the above to obtain an improved
        lower estimate on $\sG_\lambda(x)$. This procedure is repeated ad
        infinitum to obtain
        %e5.12 #&#
%e5.12 ###
\begin{equation}
                \sG_\lambda(x) \ge
                \frac{\alpha^2 \lambda}{\lambda-1}\sum_{n=0}^\infty  (
                \Upsilon(\lambda) L_\sigma^2 )^n.
        \end{equation}
        Consequently,  $\Upsilon(\lambda)\ge L_\sigma^{-2}$
        implies that $\sG_\lambda(x)=\infty$
        for all $x\in\Z^d$. If there exists a $\lambda_0>1$
        such that $\Upsilon(\lambda_0)> L_\sigma^{-2}$, then
        the preceding tells us that $\sG_{\lambda_0}\equiv\infty$.
        Now suppose, in addition, that there exists
        $z\in\Z^d$ such that
        %e5.13 #&#
%e5.13 ###
\begin{equation}\label{eqbigO}
                \E  ( | u_n(z)|^2  )=\mathrm{O}(\lambda_0^n).
        \end{equation}
        Then by the continuity of $\Upsilon$
        we can choose a finite $\lambda>\lambda_0$ such that
        $\Upsilon(\lambda) \ge L_\sigma^{-2}$, whence
        $\sG_\lambda\equiv\infty$. This yields a contradiction,
        since \eqref{eqbigO} implies that
        $\sG_\lambda(z)\le\mathrm{const}\times\sum_{n=0}^\infty (\lambda_0/\lambda)^n
        <\infty$. We have verified \eqref{eqspacetime1} when $p=2$.
        An application of H\"older's inequality proves \eqref{eqspacetime1}
        for all $p\ge2$, whence the theorem.
\end{pf*}

\begin{pf*}{Proof of Theorem \ref{thspacetimebis}}
        Because $u_0$ has finite support, it is bounded.
        Therefore, Theorem \ref{thspacetime} ensures
        the existence of an a.s.-unique solution $\mathbf u$
        to \eqref{heat}.

        Choose and fix $p\in[2 ,\infty)$,
        and let $L^p(\Z^d)$ denote the usual
        space of $p$-times summable functions $f\dvtx \Z^d\to\R$,
        normed via
        %e5.14 #&#
%e5.14 ###
\begin{equation}
                \|f\|_{L^p(\Z^d)}^p:=
                \sum_{x\in\Z^d}|f(x)|^p.
        \end{equation}
        We also define $m$ to be the counting measure
        on $\Z^d$ and consider the Banach space
        $\mathbf{B} := L^p(m\times\P)$,
        all the time noting that for all random functions $g\in\mathbf{B}$,
        %e5.15 #&#
%e5.15 ###
\begin{equation}
                \|g\|_{\mathbf{B}}
                = \biggl| \E  \biggl(\sum_{x\in\Z^d} |g(x)|^p  \biggr)
                 \biggr|^{1/p}.
        \end{equation}
        Evidently, $u_0\in L^p(\Z^d)$; we claim that
        %e5.16 #&#
%e5.16 ###
\begin{equation}\label{eqspacetimebis1}
                \frac12 \ln \Upsilon^{-1}
                 (L_\sigma^{-2} )
                \le \limsup_{n\to\infty} \frac1n  {\ln}
                 \| u_n  \|_{\mathbf{B}}
                \le \frac12 \ln \Upsilon^{-1}
                 ((c_p\lip_\sigma)^{-2} ).
        \end{equation}

        For every $\lambda>1$ we define
        $\mathbf{B}(\lambda)$ to be the Banach
        space of all $\mathcal{G}$-predictable processes~$\mathbf{f}$ with $\|\mathbf{f}\|_{\mathbf{B}(\lambda)}<\infty$,
        where
        %e5.17 #&#
%e5.17 ###
\begin{equation}
                \|\mathbf{f}\|_{\mathbf{B}(\lambda)} :=
                \sup_{n\ge 0} \lambda^{-n} \| f_n\|_{\mathbf{B}}.
        \end{equation}
        Note that $\|\mathbf{f}\|_{\lambda,p}\le\|\mathbf{f}\|_{\mathbf{B}(\lambda)}$.

        Since $u_0$ has finite support, we can use Theorem
        \ref{thfinitespt} to write
        %e5.18 #&#
%e5.18 ###
\begin{equation}
                \|u_n\|_{\mathbf{B}}
                = \mathrm{O}(n^{d/p})\times\sup_{x\in \Z^d}\|u_n(x)\|_p.
        \end{equation}
        Therefore, the following is valid for all $\lambda\in(0 ,\infty)$:
        %e5.19 #&#
%e5.19 ###
\begin{equation}
                \| \mathbf{u} \|_{\mathbf{B}(\lambda)}
                \le \mathrm{const}\cdot \sup_{n\ge 0}   \Bigl(n^{d/p}\lambda^{-n}
                \sup_{x\in\Z^d} \|u_n(x)\|_p \Bigr).
        \end{equation}
        As a result, if we select $\lambda>\lambda_0>\frac12\Upsilon^{-1}
        ((c_p\lip_\sigma)^{-2}))$, then
        %e5.20 #&#
%e5.20 ###
\begin{equation}
                \|\mathbf{u}\|_{\lambda,p}
                \le \mathrm{const}\times\|\mathbf{u}\|_{\mathbf{B}(\lambda_0)}<\infty,
        \end{equation}
        thanks to the upper bound of Theorem 2.1. It follows immediately
        from this that
        \[
        \limsup_{n\to\infty} n^{-1}{\ln}\|u_n\|_{\mathbf{B}}\le\lambda_0
        \]
        for all finite $\lambda_0>\frac12\Upsilon^{-1}
        ((c_p\lip_\sigma)^{-2})$.
        The second inequality in \eqref{eqspacetimebis1} is thus proved.
        Next we derive the first inequality in \eqref{eqspacetimebis1}.

        Thanks to Jensen's inequality, it suffices to consider only the
        case $p=2$. According to~\eqref{heat1},
        %e5.21 #&#
%e5.21 ###
\begin{equation}
                \E (  | u_{n+1}(x)  |^2 )
                \ge   | (\sP^{n+1}u_0)(x)  |^2 +
                L_\sigma^2\cdot\sum_{j=0}^n
                \sum_{y\in\Z^d}  | P^{n-j}_{x,y} |^2
                \E ( | u_j(y)  |^2 ).
        \end{equation}
        Consequently,
        %e5.22 #&#
%e5.22 ###
\begin{equation}
                \| u_{n+1}\|_{\mathbf{B}}^2
                \ge \| u_0\|_{L^2(\Z^d)}^2 + L_\sigma^2\cdot
                \sum_{j=0}^n q_{n-j}\| u_j\|_{\mathbf{B}}^2.
        \end{equation}
        We multiply both sides by $\lambda^{-(n+1)}$, then sum from
        $n=0$ to $n=\infty$ and finally apply Lemma \ref{lemUpsilon},
        in order to obtain the following:
        %e5.23 #&#
%e5.23 ###
\begin{equation}
                \sum_{n=1}^\infty\lambda^{-n}
                 \| u_n  \|_{\mathbf{B}}^2
                \ge \frac{1}{\lambda-1}\cdot\|u_0\|_{L^2(\Z^d)}^2
                + L_\sigma^2\Upsilon(\lambda)
                \cdot \sum_{k=0}^\infty
                \lambda^{-k} \| u_k \|_{\mathbf{B}}^2.
        \end{equation}
        Because $u_0\not\equiv 0$, we have $(\lambda-1)^{-1}\cdot\|u_0\|_{
        L^2(\Z^d)}^2>0$, and
        this shows that $\sum_{n=1}^\infty\lambda^{-n}\|u_n\|_{\mathbf{B}}^2
        =\infty$ whenever
        $L_\sigma^2\Upsilon(\lambda)\ge 1$.
        In particular, it must follow that $\limsup_{n\to\infty}
        \rho^{-n}\|u_n\|_{\mathbf{B}}^2=\infty$
        whenever $\rho\in(1 ,\lambda]$.
        This implies the first inequality in \eqref{eqspacetimebis1}.

        Now  we can conclude the proof from \eqref{eqspacetimebis1}.
        According to Theorem \ref{thfinitespt},
        %e5.24 #&#
%e5.24 ###
\begin{equation}
                \sup_{x\in\Z^d}\|u_n(x)\|_2\le \|u_n\|_{\mathbf{B}} \le
                \mathrm{O}(n^{d/2}) \times \sup_{x\in\Z^d}
                \|M_n\|_2.
        \end{equation}
        Therefore, \eqref{eqspacetimebis1} implies
        the theorem.
\end{pf*}

%       We record the fact that \eqref{eqspacetimebis1}
%       holds without the assumption that $\sL$ is local,
%       but it does require the condition $\sigma(0)=0$.
%       \qed

Before we prove Corollary \ref{corinter} we state and
prove an elementary convexity lemma that is due essentially to Carmona
and Molchanov \cite{CM94}, Theorem III.1.2, page 55.

%le5.3 #&#
\begin{lemma}\label{lemCM}
        Suppose $u_n(x)\ge 0$ for all $n\ge 0$, and $x\in\Z^d$,
        $\bar\gamma(p)<\infty$
        for all $p<\infty$ and $\bar\gamma(2)>0$. Then
        $\mathbf{u}$ is weakly intermittent.\vadjust{\goodbreak}
\end{lemma}

\begin{pf}
        Because $\mathbf{u}$ is non-negative,
        %e5.25 #&#
%e5.25 ###
\begin{equation}\label{eqposgam}
                \bar\gamma(\alpha) = \limsup_{n\to\infty}
                \frac1n \ln \E [ u_n(x)^\alpha ]
                \qquad\mbox{for all $\alpha\ge 0$}.\vspace*{-2pt}
        \end{equation}
        Thanks to Proposition \ref{prDuhamel},
        $\E[u_n(x)]=(\sP^n u_0)(x)$ is bounded above uniformly
        by $\sup_x u_0(x)$, which is finite. Consequently,
        %e5.26 #&#
%e5.26 ###
\begin{equation}\label{eq12}
                \bar\gamma(1)=0 < \bar\gamma(2).\vspace*{-2pt}
        \end{equation}
        Next we claim that $\bar\gamma$ is convex on $\R_+$. Indeed,
        for all $a,b\ge 0$, and $\lambda\in(0 ,1)$, H\"older's inequality
        yields the following: For all $s\in(1 ,\infty)$ with $t:=s/(s-1)$,
        %e5.27 #&#
%e5.27 ###
\begin{equation}
                \E \bigl[ u_n(x)^{\lambda a+(1-\lambda)b} \bigr] \le
                 \{ \E [ u_n(x)^{s\lambda a} ]  \}^{1/s}
                 \bigl\{ \E \bigl[ u_n(x)^{t(1-\lambda)b} \bigr]
                 \bigr\}^{1/t}.\vspace*{-2pt}
        \end{equation}
        Choose $s:=1/\lambda$
        to deduce the convexity of $\bar\gamma$ from \eqref{eqposgam}.

        Now we complete the proof: By \eqref{eq12} and convexity,
        $\bar\gamma(p)>0$ for all $p\ge 2$. If $p'>p\ge 2$, then
        we write $p=\lambda p'+(1-\lambda)$ -- with
        $\lambda := (p-1)/(p'-1)$ -- and apply convexity to conclude that
        %e5.28 #&#
%e5.28 ###
\begin{equation}\label{eqnearest}
                \bar\gamma(p) \le \lambda\bar\gamma(p') + (1-\lambda)\bar\gamma(1)
                = \frac{p-1}{p'-1}  \bar\gamma(p').\vspace*{-2pt}
        \end{equation}
        Since \eqref{eqnearest} holds, in particular, with $p\equiv 2$,
         it implies that $\bar\gamma(p')>0$.
        And the lemma follows from \eqref{eqnearest} and
        the inequality $(p-1)/(p'-1)<p/p'$.\vspace*{-3pt}
\end{pf}

\begin{pf*}{Proof of Corollary \ref{corinter}}
        Condition \ref{eqcomparison} and Theorem
        \ref{thcomparison} imply that $u_n(x)\ge 0$, and
        hence \eqref{eqposgam} holds.
        Now ``$\bar{\gamma}(2)>0$'' and ``$\bar{\gamma}(p)<\infty$
        for $p>2$'' both follow from Theorem~\ref{thspacetimebis},
        and Lemma~\ref{lemCM} completes the proof.\vspace*{-3pt}
\end{pf*}

\begin{pf*}{Proof of Theorem \ref{thtemporal}}
        The assertion about the existence and uniqueness
        of the solution to the Anderson model \eqref{heat}
        with $\sigma(z):=z$ follows from Lemma \ref{lemfirst}.
        The solution is non-negative by Lemma \ref{corcomparison}.
        Now we prove the claims about the growth of the solution~$u$.

        It is possible to check that $U_n := \sum_{x\in\Z^d} u_n(x)$
        can be written out explicitly as $U_n=  U_0\times\prod_{j=1}^n(1+\xi_j)$.
        Since $0<U_0<\infty$,
        Kolmogorov's strong law of large numbers implies that almost surely,
        %e5.29 #&#
%e5.29 ###
\begin{equation}\label{eqtemporalgamma0}
                \lim_{n\to\infty} \frac1n\ln U_n =
                \lim_{n\to\infty}\frac1n\E [\ln U_n  ]
                = \E [\ln(1+\xi_1) ]
                =\Gamma'(0^+).\vspace*{-2pt}
        \end{equation}
        Also, $\lim_{n\to\infty} n^{-1} \ln \E ( U_n^p )
        =\Gamma(p)$ for all $p\ge 0$.
        Because $ M_n\le U_n$, we have
      %e5.30 #&#
%e5.30 ###
\begin{eqnarray}\label{eqtemporalgamma0u}
                 \limsup_{n\to\infty}
                        \frac1n \ln  M_n &\le&\Gamma'(0^+)\qquad
                        \mbox{a.s.,}\nonumber\\[-2pt]
                        \limsup_{n\to\infty}
                        \frac1n \E ( \ln  M_n )
                        &\le& \Gamma'(0^+)  \quad  \mbox{and}\\[-2pt]
                 \limsup_{n\to\infty} \frac1n\ln\E (
                        M_n^p )&\le&\Gamma(p)
                        \qquad\mbox{for all $p\in[0 ,\infty)$}.
        \nonumber\vspace*{-2pt}
\end{eqnarray}
        Next we strive to establish the complementary inequalities
        to these.\vadjust{\goodbreak}

        In order to derive the second, and final,
        half of the theorem, we choose and fix
        some $x_0\in\Z^d$ such that $u_0(x_0)>0$.
        Let $\mathbf{v}:=\{v_n(x)\}_{n\ge 0,x\in\Z^d}$
        solve \eqref{heat} with $\sigma(z)=z$,
        subject to $v_0(x)=u_0(x_0)$
        if $x=x_0$ and $v_0(x)=0$ otherwise.
        The existence and uniqueness of $v$ follows from Lemma \ref{lemfirst}.
        By Corollary \ref{corcomparison},
        %e5.31 #&#
%e5.31 ###
\begin{equation}\label{equv}
                0\le v_n(x)\le u_n(x)\qquad\mbox{for all $n\ge 0$
                and $x\in\Z^d$.}
        \end{equation}
        Let $V_n:=\sum_{x\in\Z^d} v_n(x)$.
        Then
        %e5.32 #&#
%e5.32 ###
\begin{equation}\label{eqtemporalgamma0v}
                \lim_{n\to\infty} \frac1n \ln V_n =
                \lim_{n\to\infty} \frac1n \E [\ln V_n  ]
                = \E [\ln(1+\xi_1) ] =\Gamma'(0^+).
        \end{equation}
        Also, $\lim_{n\to\infty} n^{-1} \ln \E ( V_n^p )
        =\Gamma(p)$ for all $p\in[0 ,\infty)$.
        Recall $R$ from \eqref{defR}. Because $v_0=0$ off of $\{x_0\}$,
        \eqref{eqgrowspt} implies that
        $V_n \le 2^d\{nR+1\}^d\times \sup_{x\in\Z^d} v_n(x)$.
        Owing to Theorem \ref{thcomparison}, $\sup_{x\in\Z^d}
        v_n(x)\le M_n$. Therefore,
  %e5.33 #&#
%e5.33 ###
\begin{eqnarray}
                \liminf_{n\to\infty} \frac1n \ln  M_n &\ge& \Gamma'(0^+),\qquad
                        \liminf_{n\to\infty} \frac1n\E ( \ln  M_n  )\ge
                        \Gamma'(0^+)\quad\mbox{and}\nonumber
                        \\[-8pt]
                        \\[-8pt]
               \liminf_{n\to\infty} \frac1n \ln
                        \E (  M_n^p )
                        &\ge&\Gamma(p)\qquad\mbox{for all $p\in[0 ,\infty)$}.
        \nonumber
\end{eqnarray}
        Together with \eqref{eqtemporalgamma0u},
        these bounds prove Theorem \ref{thtemporal}.
\end{pf*}

%s6 #&#
%s6 ###
\section{An example}\label{sec6}

Let us consider \eqref{heat} in the special case
that: (i) $\xi$s are independent mean-zero
variance-one random variables; (ii) $\sigma(z)=\nu z$
for a fixed $\nu >0$; (iii) $u_0$ has finite support;
and (iv)~$\sL$ is the generator of a simple symmetric random
walk on $\Z$.
That is,
%e6.1 #&#
%e6.1 ###
\begin{equation}
        (\sL h)(x) = \frac{h(x+1)+h(x-1)-2h(x)}{2},
\end{equation}
for every function $h\dvtx \Z\to\R$ and all $x\in\Z$.
The operator
$2\sL$ is called the \emph{graph Laplacian}
on $\Z$, and the resulting form,
%e6.2 #&#
%e6.2 ###
\begin{equation}
        u_{n+1}(x)-u_n(x) = (\sL u_n)(x) + \nu u_n(x)\xi_n(x),
\end{equation}
of \eqref{heat} is an \emph{Anderson model} of a parabolic type
\cite{CM94,Molch91}.
Theorems \ref{thspacetime} and \ref{thspacetimebis}
together imply that the upper Lyapounov exponent
of the solution to \eqref{heat} is $\ln\Upsilon^{-1} (\nu^{-2} )$
in this case.
We compute the quantity $\Upsilon^{-1}(\nu^{-2})$ next.
The following might suggest that one cannot hope to compute upper
Lyapounov exponents explicitly in general.

%pr6.1 #&#
\begin{proposition}\label{prU}
        If $\nu>0$, then
        %e6.3 #&#
%e6.3 ###
\begin{equation}
                \Upsilon^{-1} ( \nu^{-2} )=\inf
                \biggl \{ \lambda>1\dvt
                {}_1F_1 \biggl (\frac12  ;1 ; \frac1\lambda \biggr)<
                \frac{\lambda}{\nu^2}  \biggr\}.
        \end{equation}
\end{proposition}

\begin{pf}
        Recall the $q_k$s from \eqref{eqq}.
        According to Plancherel's theorem and symmetry,
        %e6.4 #&#
%e6.4 ###
\begin{equation}
                q_n =\frac1\uppi \int_0^\uppi
                \biggl ( \frac{1+\cos(2\xi)}{2}
         \biggr)^n\,\d\xi.
        \end{equation}
        We may apply the half-angle formula
        for cosines, and then Wallis's formula
        (Davis \cite{PJ-Davis}, equation ({6.1.49}), page 258),
        in order to find that if $n\ge 1$, then
        %e6.5 #&#
%e6.5 ###
\begin{equation}
                q_n = \frac{(2n-1)!!}{(2n)!!},
        \end{equation}
        where ``$!!$'' denotes the double factorial.
        Therefore, ${}_1F_1(1/2 ;1 ;\bullet)$
        is the generating function of the sequence
        $\{q_n\}_{n=0}^\infty$; confer with Slater
         \cite{Slater}, equation ({13.1.2}), page 504.
%       Consequently,
%       \begin{equation}
%               \sum_{n=0}^\infty \lambda^{-n}q_n
%                = {}_1F_1
%                ( \frac12 ;1 ; \frac{1}{\lambda}
%                ).
%       \end{equation}
        This and Lemma \ref{lemUpsilon} together prove that\vspace*{-1.3pt}
        %e6.6 #&#
%e6.6 ###
\begin{equation}
                \lambda\Upsilon(\lambda) =  {}_1F_1
                \biggl ( \frac12  ;1 ;  \frac1\lambda \biggr),
        \end{equation}
        and the lemma follows since
        $\Upsilon$ is a continuous and strictly
        decreasing function on $(0 ,\infty)$.\vspace*{-1.3pt}
\end{pf}

\section*{Acknowledgements}
An anonymous referee
has kindly made a number of detailed suggestions and corrections
which led to a much better paper. We thank him/her wholeheartedly.
Research supported in part by NSF Grant DMS-07-04024.\vspace*{-1.3pt}

% imsref loaded by smiklovaite, 2011-12-03 08:59:45
% imsref loaded by smiklovaite, 2011-12-03 09:09:10

\printhistory

\end{document}